\documentclass[10pt, a4paper]{article}

\usepackage{amsfonts, latexsym, amsmath, amssymb, fancyhdr, amsthm, graphicx, subfigure}
\usepackage[english]{babel}

\newtheorem{theorem}{Theorem}[section]
\newtheorem{lemma}{Lemma}[section]
\newtheorem{proposition}{Proposition}[section]

\newtheorem{corollary}{Corollary}[section]
\newtheorem{remark}{Remark}

\title{\textbf{Stability and Hopf Bifurcation in a \\ Mathematical Model
of \\ Pluripotent Stem Cell Dynamics}}
\author{Mostafa Adimy$^{\dagger}$, \quad Fabien Crauste$^{\dagger}$ \quad and \quad Shigui Ruan\footnote{\small Research was partially supported by the NSERC of Canada
and the College of Arts and Sciences at the University of Miami.
On leave from Dalhousie University, Halifax, Canada.} }
\date{Year 2004}

\begin{document}

\maketitle

\begin{center}
{\large $^{\dagger}$}\emph{Laboratoire de Math\'ematiques Appliqu\'ees, FRE 2570,} \\
\emph{Universit\'e de Pau et des Pays de l'Adour,}\\
\emph{Avenue de l'universit\'e, 64000 Pau, France.}\\
\emph{E-mail: mostafa.adimy@univ-pau.fr, fabien.crauste@univ-pau.fr}\\
%\emph{E-mail: fabien.crauste@univ-pau.fr}\\
\quad\\
{\large $^*$}\emph{Department of Mathematics, University of
Miami,}\\\emph{P. O. Box 249085, Coral Gables, FL 33124-4250,
USA.}\\\emph{ E-mail: ruan@math.miami.edu}
\end{center}

\quad

\begin{abstract}
We study a mathematical model describing the dynamics of a
pluripotent stem cell population involved in the blood production
process in the bone marrow. This model is a differential equation
with a time delay. The delay describes the cell cycle duration and
is uniformly distributed on an interval. We obtain stability
conditions independent of the delay. We also show that the
distributed delay can destabilize the entire system. In
particularly, it is shown that Hopf bifurcations can occur.
\end{abstract}

\bigskip{}


%\newpage
%%%%%%%%%%%%%%%%%%%%%%%%%%%%%%%%%%%%%%%%%%%%%%%%%%%%%%%%%%%%%%%%%%%%%%%%%%%%%%%%%%%%%%%%
%%%%%%%%%%%%%%%%%%%%%%%%%%%%%%%%%%%%%%%%%%%%%%%%%%%%%%%%%%%%%%%%%%%%%%%%%%%%%%%%%%%%%%%%
%%%%%%%%%%%%%%%%
%%%%%%%%%%%%%%%%                       INTRODUCTION-MOTIVATION
%%%%%%%%%%%%%%%%
%%%%%%%%%%%%%%%%%%%%%%%%%%%%%%%%%%%%%%%%%%%%%%%%%%%%%%%%%%%%%%%%%%%%%%%%%%%%%%%%%%%%%%%%
%%%%%%%%%%%%%%%%%%%%%%%%%%%%%%%%%%%%%%%%%%%%%%%%%%%%%%%%%%%%%%%%%%%%%%%%%%%%%%%%%%%%%%%%
\section{Introduction}

Blood production process, called hematopoiesis, is one of the
major biological phenomena occurring in human body. It takes place
in the bone marrow where pluripotent stem cells give birth to
mature cells. After ejecting their nuclei, these cells enter the
bloodstream and become blood cells.

According to the study of Burns and Tannock \cite{burns}, the
population of pluripotent stem cells can be divided into two
distinct groups: quiescent cells and proliferating cells.
Mathematical models describing the dynamics of this cell
population have been studied since the end of the seventies, in
particularly by Mackey \cite{mackey1978, mackey1979}. We refer to
the review articles by Haurie {\it et al.} \cite{hdm98} and Mackey
{\it et al.} \cite{mackey2003} for further study and more
references on this topic. More recently, Pujo-Menjouet et al.
\cite{mackeypujo2} and Pujo-Menjouet and Mackey \cite{mackeypujo}
proved the existence of a Hopf bifurcation for the hematopoiesis
model proposed in \cite{mackey1978}. In all these works, the
authors assumed that the proliferating phase duration is constant.
Mathematically, this means that the delay in their models is a
discrete delay. However, experimental data (see Bradford \emph{et
al.} \cite{bradford}) indicate that cells do not spend the same
time in the proliferating phase.

In this paper, taking into account this assumption, we assume that
the delay (or proliferating phase duration) is uniformly
distributed on an interval. The main objective is to investigate
the effect of time delay on the dynamical solutions. It is shown
that there exist some critical values of time delay such that a
local Hopf bifurcation occurs at the non-trivial equilibrium.

The paper is organized as follows. In section \ref{sectionmodel},
we present our model, which is given in equation
(\ref{equationx}). In section \ref{sectiongammazero}, we derive
stability conditions for the two equilibria of equation
(\ref{equationx}) which do not depend on the delay. We show the
existence of Hopf bifurcations at the non-trivial equilibrium in
section \ref{hopfbifurcation}. A brief discussion is given in
section \ref{discussion}.

%%%%%%%%%%%%%%%%%%%%%%%%%%%%%%%%%%%%%%%%%%%%%%%%%%%%%%%%%%%%%%%%%%%%%%%%%%%%%%%%%%%%%%%%
%%%%%%%%%%%%%%%%%%%%%%%%%%%%%%%%%%%%%%%%%%%%%%%%%%%%%%%%%%%%%%%%%%%%%%%%%%%%%%%%%%%%%%%%
%%%%%%%%%%%%%%%%
%%%%%%%%%%%%%%%%                          MATHEMATICAL MODEL
%%%%%%%%%%%%%%%%
%%%%%%%%%%%%%%%%%%%%%%%%%%%%%%%%%%%%%%%%%%%%%%%%%%%%%%%%%%%%%%%%%%%%%%%%%%%%%%%%%%%%%%%%
%%%%%%%%%%%%%%%%%%%%%%%%%%%%%%%%%%%%%%%%%%%%%%%%%%%%%%%%%%%%%%%%%%%%%%%%%%%%%%%%%%%%%%%%
\section{The Model}\label{sectionmodel}

Pluripotent stem cells can be either in a resting phase, also
known as $G_0$-phase, or in a proliferating phase. In the resting
phase, they can die at a constant rate $\delta\geq0$, which also
includes the cellular differentiation, or be introduced in the
proliferating phase at a rate $\beta$. According to the work of
Sachs \cite{sachs}, $\beta$ is assumed to depend on the resting
phase population.

In the proliferating phase, which is in fact the so-called cell
cycle, pluripotent stem cells are committed to divide and give
birth to two daughter cells at the end of this phase. The two
daughter cells enter directly the resting phase and complete the
cycle. We assume that proliferating cells divide according to a
uniform law $f$ on an interval $[\tau_{min},\tau]$ with
$0\leq\tau_{min}<\tau<+\infty$. This assumption comes from the
fact that, even if only a little is known about phenomena involved
in hematopoiesis, there are strong evidences (see Bradford
\emph{et al.} \cite{bradford}) indicating that cells do not divide
at the same age. The function $f$ is then defined by
\begin{displaymath}
f(r)=\left\{\begin{array}{ll}
\displaystyle\frac{1}{\tau-\tau_{min}},&\qquad \textrm{ if }
r\in[\tau_{min},\tau],\\
0,& \qquad \textrm{ otherwise.}
\end{array}\right.
\end{displaymath}

Let $x(t)$ denote the pluripotent stem cell population density
(cells/kg) at time $t\geq0$. It satisfies the nonlinear delay
differential equation
\begin{equation} \label{equationx}
x^{\prime}(t)=-\big( \delta+\beta(x(t)) \big)x(t) +
\frac{2}{\tau-\tau_{min}}\int_{\tau_{min}}^{\tau}
\beta(x(t-r))x(t-r) dr.
\end{equation}
The first term in the right-hand side of equation
(\ref{equationx}) accounts for the cellular loss due to mortality
and cellular differentiation, $\delta x(t)$, and introduction in
the cell cycle, $\beta(x(t))x(t)$. The second term is for the
division of proliferating cells into two daughter cells during
mitosis. Proliferating cells are in fact resting cells introduced
in the proliferating phase one generation earlier, so that the
quantity $\beta(x(t-r))x(t-r)$ appears with a time delay. The
factor 2 is, of course, for the division of each proliferating
cell into two daughter cells.

In the following, the rate of reintroduction in the proliferating
compartment $\beta$ is taken to be a monotone and decreasing Hill
function, given by
\begin{displaymath}
\beta(x)=\beta_0\frac{\theta^n}{\theta^n+x^n} \qquad \textrm{ for
} x\geq0.
\end{displaymath}
The coefficient $\beta_0>0$ is the maximum rate of reintroduction,
$\theta\geq0$ is the $G_0$-phase population density for which the
rate of re-entry $\beta$ attains its maximum rate of change with
respect to the resting phase population, and $n\geq0$ describes
the sensitivity of $\beta$ with changes in the population. This
function was firstly used in hematopoiesis models by Mackey
\cite{mackey1978} in 1978.

In \cite{mackey1978} and \cite{mackey2001}, Mackey gave values of
the above parameters for a normal human body production. These
values are
\begin{equation}\label{parametersvalues}
\delta=0.05\textrm{ d}^{-1}, \quad \beta_0=1.77\textrm{ d}^{-1}
\quad \textrm{ and } \quad n=3.
\end{equation}
The value of $\theta$ is usually $\theta=1.62\times 10^8\textrm{
cells/kg}.$ However, since we shall study the qualitative behavior
of the pluripotent stem cells population, the value of $\theta$ is
not really important and could be normalized without loss of
generality.

Now if we consider an initial continuous nonnegative function
$\varphi$ defined on $[-\tau,0]$, then the equation
(\ref{equationx}) has a unique continuous and nonnegative solution
$x^{\varphi}(t)$, defined for $t\geq -\tau$, such that
\begin{displaymath}
x^{\varphi}(s)=\varphi(s) \qquad \textrm{ for } s\in[-\tau,0].
\end{displaymath}
This can be obtained by using the results in Hale and Verduyn
Lunel \cite{halelunel}.

Notice that equation (\ref{equationx}) has at most two equilibria,
the trivial equilibrium $x\equiv 0$ and a non-trivial positive
equilibrium $x\equiv x^*$. The trivial equilibrium always exists
and corresponds to the extinction of the population.

%%
%%      EXISTENCE OF THE NON-TRIVIAL EQUILIBRIUM
%%
\begin{proposition}\label{existxstar}
Equation (\ref{equationx}) has a non-trivial positive equilibrium
$x\equiv x^*$ if and only if
\begin{equation}\label{existencexstar}
\beta_0>\delta>0.
\end{equation}
In this case, $x^*$ is explicitly given by
\begin{displaymath}
x^*=\theta\bigg(\frac{\beta_0}{\delta}-1\bigg)^{1/n}.
\end{displaymath}
\end{proposition}

\begin{proof}
Let $x^*$ be an equilibrium of equation (\ref{equationx}). Then
$x^*$ satisfies
\begin{displaymath}
x^*\big(\beta(x^*)-\delta\big)=0.
\end{displaymath}
Consequently, equation (\ref{equationx}) has a non-trivial
equilibrium if and only if the equation
\begin{displaymath}
\beta(x^*)=\delta
\end{displaymath}
has a  non-trivial solution. Since the function $\beta$ is
decreasing and positive with $\beta(0)=\beta_0$, then equation
(\ref{equationx}) has a non-trivial equilibrium if and only if
condition (\ref{existencexstar}) holds.
\end{proof}

In the next section, we shall study the stability of the two
equilibria of equation (\ref{equationx}).

%%%%%%%%%%%%%%%%%%%%%%%%%%%%%%%%%%%%%%%%%%%%%%%%%%%%%%%%%%%%%%%%%%%%%%%%%%%%%%%%%%%%%%%%
%%%%%%%%%%%%%%%%%%%%%%%%%%%%%%%%%%%%%%%%%%%%%%%%%%%%%%%%%%%%%%%%%%%%%%%%%%%%%%%%%%%%%%%%
%%%%%%%%%%%%%%%%%
%%%%%%%%%%%%%%%%%                LOCAL STABILITY
%%%%%%%%%%%%%%%%%
%%%%%%%%%%%%%%%%%%%%%%%%%%%%%%%%%%%%%%%%%%%%%%%%%%%%%%%%%%%%%%%%%%%%%%%%%%%%%%%%%%%%%%%%
%%%%%%%%%%%%%%%%%%%%%%%%%%%%%%%%%%%%%%%%%%%%%%%%%%%%%%%%%%%%%%%%%%%%%%%%%%%%%%%%%%%%%%%%
\section{Stability}\label{sectiongammazero}

Throughout this section, we are interested in the stability of the
equilibria of equation (\ref{equationx}), in particularly the
stability of the non-trivial equilibrium $x\equiv x^*$. We start
by giving a result on the global stability of the trivial
equilibrium of (\ref{equationx}).

%%
%%      GLOBAL STABILITY OF THE TRIVIAL EQUILIBRIUM
%%
\begin{theorem}\label{stabilitytrivialequilibrium}
The trivial equilibrium $x\equiv0$ of equation (\ref{equationx})
is globally stable if
\begin{displaymath}
\beta_0<\delta.
\end{displaymath}
\end{theorem}

\begin{proof}
The proof uses a similar technique employed by Adimy and Crauste
\cite{adimycrauste2003}. It is based on the construction of a
Lyapunov functional.

Denote by $C^+$ the space of all continuous nonnegative functions
on $[-\tau,0]$. Let $B$ be the function defined by
\begin{displaymath}
B(x)=\int_0^x \beta(s)s \ ds \qquad \textrm{ for } x\geq0.
\end{displaymath}
Consider the mapping $J: C^+\to[0,+\infty)$ defined, for
$\varphi\in C^+$, by
\begin{displaymath}
J(\varphi)=B(\varphi(0))+\frac{1}{\tau-\tau_{min}}\int_{\tau_{min}}^{\tau}
\int_{-r}^{0}\Big(\beta(\varphi(a))\varphi(a)\Big)^2 da dr.
\end{displaymath}
Then,
\begin{displaymath}
\dot J(\varphi)=\dot \varphi(0)\beta(\varphi(0)) \varphi(0)
+\frac{1}{\tau-\tau_{min}}
\int_{\tau_{min}}^{\tau}\big(\beta(\varphi(0))\varphi(0)\big)^{2}
-\big(\beta(\varphi(-r))\varphi(-r)\big)^{2}dr.
\end{displaymath}
Since
\begin{displaymath}
\dot \varphi(0) = -\big(\delta+\beta(\varphi(0))\big)\varphi(0) +
\frac{2}{\tau-\tau_{min}}\int_{\tau_{min}}^{\tau}\beta(\varphi(-r))\varphi(-r)dr,
\end{displaymath}
we obtain that
\begin{displaymath}
\begin{array}{rcl}
\dot J(\varphi) & = &
-\big(\delta+\beta(\varphi(0))\big)\beta(\varphi(0))\varphi^{2}(0)
+\displaystyle\frac{2}{\tau-\tau_{min}}\int_{\tau_{min}}^{\tau}\big(\beta(\varphi(0))\varphi(0)\big)^2dr \vspace{1ex}\\
        &   & -\displaystyle\frac{1}{\tau-\tau_{min}}\int_{\tau_{min}}^{\tau}\big(\beta(\varphi(0))\varphi(0)-\beta(\varphi(-r))\varphi(-r)\big)^{2}dr.
\end{array}
\end{displaymath}
Hence,
\begin{displaymath}
\dot J(\varphi) \leq
-\big(\delta-\beta(\varphi(0))\big)\beta(\varphi(0))\varphi(0)^2.
\end{displaymath}
Let $\alpha$ be the function defined, for $x \geq 0$, by
\begin{displaymath}
\alpha(x)=(\delta-\beta(x))\beta(x)x^{2}.
\end{displaymath}
Assume that $\beta_0<\delta.$ Since $\beta$ is a decreasing
function, it follows that the function $x\mapsto \delta-\beta(x)$
is positive for $x\geq0$. Hence, $\alpha$ is nonnegative on
$[0,+\infty)$ and $\alpha(x)=0$ if and only if $x=0$.
Consequently, the mapping $J$ is a Lyapunov functional when
$\beta_0<\delta$. We then deduce that the trivial equilibrium of
(\ref{equationx}) is globally stable.
\end{proof}

The result in Theorem \ref{stabilitytrivialequilibrium} describes
the fact that when $x\equiv 0$ is the only equilibrium of
(\ref{equationx}), the population is doomed to extinction except
when $\beta_0=\delta$.

Now we focus on the stability of the positive equilibrium $x\equiv
x^*$ of equation (\ref{equationx}). To ensure the existence of the
equilibrium $x\equiv x^*$, we assume that condition
(\ref{existencexstar}) holds; that is,
\begin{displaymath}
\beta_0>\delta>0.
\end{displaymath}
We do not expect to obtain conditions for the global stability of
$x\equiv x^*$. However, local stability results can be obtained by
linearizing equation (\ref{equationx}) about $x^*$. Set
\begin{equation}\label{betastargammazero}
\beta^*:=\frac{d}{dx}\Big(\beta(x)x\Big)\Big|_{x=x^*}
=\delta\bigg(1-n\frac{\beta_0-\delta}{\beta_0}\bigg).
\end{equation}
The linearization of equation (\ref{equationx}) at $x^*$ is
\begin{displaymath}
x^{\prime}(t)=-(\delta+\beta^*)x(t)+\frac{2\beta^*}{\tau-\tau_{min}}
\int_{\tau_{min}}^{\tau}x(t-r)dr.
\end{displaymath}
The characteristic equation of (\ref{equationx}) is given by
\begin{equation}\label{characteristicequationzero}
\Delta(\lambda):=\lambda+\delta+\beta^*-\frac{2\beta^*}
{\tau-\tau_{min}}\int_{\tau_{min}}^{\tau}e^{-\lambda r}dr=0.
\end{equation}

We now state and prove our first result on the stability of
$x\equiv x^*$.

%%
%%      THEOREM: STABILITY OF x* WHEN \beta^*\geq0
%%
\begin{theorem}\label{stabilitygammazero}
Assume that
\begin{displaymath}
n\frac{\beta_0-\delta}{\beta_0}\leq1.
\end{displaymath}
Then the non-trivial equilibrium $x\equiv x^*$ of equation
(\ref{equationx}) is locally asymptotically stable.
\end{theorem}

\begin{proof}
We show, in fact, that $x^*$ is stable when $\beta^*\geq0$. By the
definition of $\beta^*$ given by (\ref{betastargammazero}), it
follows that
\begin{displaymath}
\beta^*\geq0 \qquad \textrm{ if and only if } \qquad
n\frac{\beta_0-\delta}{\beta_0}\leq1.
\end{displaymath}
So we assume that $\beta^*\geq0$.

We first assume that $\Delta(\lambda)$, given by
(\ref{characteristicequationzero}), is a real function. Then,
$\Delta(\lambda)$ is continuously differentiable and its first
derivative is given by
\begin{equation}\label{firstderivativezero}
\frac{d\Delta}{d\lambda}(\lambda)=1+\frac{2\beta^*}{\tau-\tau_{min}}
\int_{\tau_{min}}^{\tau} re^{-\lambda r}dr.
\end{equation}
One can see that $d\Delta/d\lambda$ is positive for
$\lambda\in\mathbb{R}$ as soon as $\beta^*\geq 0$. Moreover,
\begin{displaymath}
\lim_{\lambda\to -\infty} \Delta(\lambda)=-\infty \qquad \textrm{
and } \qquad \lim_{\lambda\to +\infty} \Delta(\lambda)=+\infty.
\end{displaymath}
Consequently, $\Delta(\lambda)$ has a unique real root
$\lambda_0$. Since
\begin{displaymath}
\Delta(0)=\delta-\beta^*=n\delta\Big(1-\frac{\delta}{\beta_0}\Big)>0,
\end{displaymath}
we deduce that $\lambda_0<0$.

Now, we show that if $\lambda$ is a characteristic root of
equation (\ref{characteristicequationzero}), then
$Re(\lambda)\leq\lambda_0$. By contradiction, we assume that there
exists a characteristic root $\lambda=\mu+i\omega$ of equation
(\ref{characteristicequationzero}) such that $\mu>\lambda_0$. By
considering the real part of $\Delta(\lambda)$, we obtain that
\begin{displaymath}
\mu+\delta+\beta^*-\frac{2\beta^*}{\tau-\tau_{min}}\int_{\tau_{min}}^{\tau}e^{-\mu
r}\cos(\omega r) dr=0.
\end{displaymath}
Consequently,
\begin{displaymath}
\mu-\lambda_0=\frac{2\beta^*}{\tau-\tau_{min}}\int_{\tau_{min}}^{\tau}\big(e^{-\mu
r}\cos(\omega r)-e^{-\lambda_0 r }\big)dr<0.
\end{displaymath}
This yields a contradiction. We conclude that every characteristic
root $\lambda$ of (\ref{characteristicequationzero}) is such that
$Re(\lambda)\leq\lambda_0$. Hence, all characteristic roots of
(\ref{characteristicequationzero}) have negative real parts and
the equilibrium $x\equiv x^*$ is locally asymptotically stable.
\end{proof}

When $\beta^*<0$, that is, when
\begin{displaymath}
1<n\frac{\beta_0-\delta}{\beta_0},
\end{displaymath}
the stability cannot occur for all values of $\tau_{min}$ and
$\tau$. In particularly, we shall show that a Hopf bifurcation can
occur (see Theorem \ref{theoremhopfbifurcation}). However, we can
still have the stability of the non-trivial equilibrium $x\equiv
x^*$ for values of $n$, $\beta_0$ and $\delta$ if
$n(\beta_0-\delta)/\beta_0$ is not too large. This will be
considered in the next theorem.

To present the results, without loss of generality we assume that
\begin{displaymath}
\tau_{min}=0.
\end{displaymath}
We want to point out that the results we are going to show remain
true when $\tau_{min}>0$, but the proof is more complicated.

Define a function $K$, for $x\geq0$, by
\begin{equation}\label{functionK}
K(x)=\frac{\sin(x)}{x}
\end{equation}
and let $x_1$ be the unique solution of the equation
\begin{displaymath}
x_1=\tan(x_1),\qquad x_1\in(\pi,\frac{3\pi}{2}).
\end{displaymath}
Set
\begin{displaymath}
u_0:=\cos(x_1) \in (-1,0).
\end{displaymath}
Then $K^{\prime}(x_1)=0$ and
\begin{equation}\label{hypmin}
u_0=K(x_1)=\min_{x\geq0}K(x).
\end{equation}

We have the following local stability theorem.

%%
%%      THEOREM: STABILITY OF x* WHEN \beta^*<0
%%
\begin{theorem}\label{theoremstab2}
Assume that
\begin{equation}\label{hypbeta}
1<n\frac{\beta_0-\delta}{\beta_0}<\frac{2(1-u_0)}{1-2u_0}.
\end{equation}
Then the non-trivial equilibrium $x\equiv x^*$ of equation
(\ref{equationx}) is locally asymptotically stable.
\end{theorem}

\begin{proof}
Let us assume that (\ref{hypbeta}) holds. Then $\beta^*<0$,
$\delta+\beta^*>0$ and
\begin{equation}\label{hypu0}
\frac{\delta+\beta^*}{2\beta^*}<u_0.
\end{equation}

\noindent By contradiction, assume that there exists a
characteristic root $\lambda=\mu+i\omega$ of
(\ref{characteristicequationzero}) with $\mu>0$. Then,
\begin{displaymath}
\mu=-(\delta+\beta^*)+\frac{2\beta^*}{\omega\tau}\int_0^{\omega\tau}
e^{-\displaystyle\frac{\mu}{\omega}r}\cos(r) dr
\end{displaymath}
and
\begin{displaymath}
\omega=-\frac{2\beta^*}{\omega\tau}\int_0^{\omega\tau}
e^{-\displaystyle\frac{\mu}{\omega}r}\sin(r) dr.
\end{displaymath}
Integrating by parts, we obtain that
\begin{displaymath}
2\mu=-(\delta+\beta^*)+2\beta^*e^{-\mu\tau}K(\omega\tau).
\end{displaymath}
Consequently,
\begin{displaymath}
\mu<-(\delta+\beta^*)+2\beta^*e^{-\mu\tau}K(\omega\tau).
\end{displaymath}
If $\omega\tau$ is such that
\begin{displaymath}
\sin(\omega\tau)\geq0,
\end{displaymath}
then
\begin{displaymath}
2\beta^*K(\omega\tau)<0
\end{displaymath}
and
\begin{displaymath}
\mu<-(\delta+\beta^*)\leq0.
\end{displaymath}
So we obtain a contradiction.

Similarly, if $\omega\tau$ is such that
\begin{displaymath}
\sin(\omega\tau)<0,
\end{displaymath}
then, from (\ref{hypmin}) and (\ref{hypu0}), we deduce that
\begin{displaymath}
\frac{\delta+\beta^*}{2\beta^*}\leq K(\omega\tau).
\end{displaymath}
It implies that
\begin{displaymath}
\delta+\beta^*\geq
2\beta^*K(\omega\tau)>2\beta^*e^{-\mu\tau}K(\omega\tau).
\end{displaymath}
Therefore,
\begin{displaymath}
\mu<-(\delta+\beta^*)+2\beta^*K(\omega\tau)\leq0.
\end{displaymath}
Again we obtain a contradiction. Hence, all characteristic roots
$\lambda$ of (\ref{characteristicequationzero}) are such that
$Re(\lambda)\leq0$.

Now, we assume that (\ref{characteristicequationzero}) has a
purely imaginary characteristic root $\lambda=i\omega$. Then
$\omega$ and $\tau$ satisfy
\begin{equation}\label{eqpureimag}
K(\omega\tau)=\frac{\delta+\beta^*}{2\beta^*}.
\end{equation}
Using (\ref{hypmin}) and (\ref{hypu0}), we obtain a contradiction.
Consequently, (\ref{eqpureimag}) has no solution and equation
(\ref{characteristicequationzero}) does not have purely imaginary
roots. We conclude that all characteristic roots of
(\ref{characteristicequationzero}) have negative real parts and
$x\equiv x^*$ is locally asymptotically stable.
\end{proof}

From Theorems \ref{stabilitygammazero} and \ref{theoremstab2}, it
follows that the non-trivial equilibrium $x\equiv x^*$ of equation
(\ref{equationx}) is locally asymptotically stable when
\begin{equation}\label{stabilityn}
0\leq n\frac{\beta_0-\delta}{\beta_0}<\frac{2(1-u_0)}{1-2u_0}.
\end{equation}
We are going to show that as soon as condition (\ref{stabilityn})
does not hold, then the equilibrium can be destabilized. In the
next section, we shall show that if condition (\ref{hypu0}) does
not hold, then a Hopf bifurcation indeed occurs at $x\equiv x^*$.

%%%%%%%%%%%%%%%%%%%%%%%%%%%%%%%%%%%%%%%%%%%%%%%%%%%%%%%%%%%%%%%%%%%%%%%%%%%%%%%%%%%%%%%%
%%%%%%%%%%%%%%%%%%%%%%%%%%%%%%%%%%%%%%%%%%%%%%%%%%%%%%%%%%%%%%%%%%%%%%%%%%%%%%%%%%%%%%%%
%%%%%%%%%%%%%%%%%
%%%%%%%%%%%%%%%%%                HOPF BIFURCATIONS
%%%%%%%%%%%%%%%%%
%%%%%%%%%%%%%%%%%%%%%%%%%%%%%%%%%%%%%%%%%%%%%%%%%%%%%%%%%%%%%%%%%%%%%%%%%%%%%%%%%%%%%%%%
%%%%%%%%%%%%%%%%%%%%%%%%%%%%%%%%%%%%%%%%%%%%%%%%%%%%%%%%%%%%%%%%%%%%%%%%%%%%%%%%%%%%%%%%
\section{Hopf Bifurcations}\label{hopfbifurcation}

In this section we are going to show that the non-trivial
equilibrium $x\equiv x^*$ of equation (\ref{equationx}) can be
destabilized \emph{via} Hopf bifurcations. The time delay $\tau$
will be used as a bifurcation parameter. This result is obtained
in Theorem \ref{theoremhopfbifurcation}.

Recall that the non-trivial equilibrium $x\equiv x^*$ of equation
(\ref{equationx}) exists if and only if $\beta_0>\delta>0.$ In the
following, without loss of generality we assume that
\begin{displaymath}
\tau_{min}=0.
\end{displaymath}
Again the results still hold  when $\tau_{min}>0$, but the proof
is easier to understand when $\tau_{min}=0$.

We look for purely imaginary roots of $\Delta(\lambda)$. Of
course, we assume that $\beta^*<0$, otherwise $x\equiv x^*$ is
locally asymptotically stable. Let $\lambda=i\omega$, with
$\omega\in\mathbb{R}$, be a purely imaginary characteristic root
of equation (\ref{characteristicequationzero}). Then, $\tau$ and
$\omega$ satisfy the following system
\begin{equation}\label{system}
\left\{\begin{array}{rcl}
\delta+\beta^*\big(1-2C(\tau,\omega)\big)&=&0,\\
\omega+2\beta^*S(\tau,\omega)&=&0,
\end{array}\right.
\end{equation}
where
\begin{displaymath}
C(\tau,\omega)=\frac{1}{\tau}\int_{0}^{\tau}\cos(\omega
r)dr,\qquad
S(\tau,\omega)=\frac{1}{\tau}\int_{0}^{\tau}\sin(\omega r)dr.
\end{displaymath}
First, one can see that $\omega=0$ cannot be a solution of
(\ref{system}). Otherwise
\begin{displaymath}
\delta=\beta^*<0.
\end{displaymath}
Moreover, if $\omega$ is a solution of system (\ref{system}), then
$-\omega$ is also a solution of (\ref{system}). Hence, we only
look for positive solutions $\omega$.

One can check that $C(\tau,\omega)$ and $S(\tau,\omega)$ are
given, for $\tau>0$ and $\omega>0$, by
\begin{displaymath}
C(\tau,\omega)=\frac{\sin(\omega\tau)}{\omega\tau}=K(\omega\tau),
\qquad S(\tau,\omega)=\frac{1-\cos(\omega\tau)}{\omega\tau},
\end{displaymath}
where the function $K$ is defined by (\ref{functionK}).
Consequently, system (\ref{system}) can be rewritten as
\begin{eqnarray}
K(\omega\tau)&=&\frac{\delta+\beta^*}{2\beta^*},\vspace{1ex}\label{eqsin}\\
\displaystyle\frac{\cos(\omega\tau)-1}{(\omega\tau)^2}&=&\frac{1}{2\beta^*\tau}.\label{eqcos}
\end{eqnarray}

\begin{figure}[!hpt]
\begin{center}
\includegraphics[width=10cm,height=8cm]{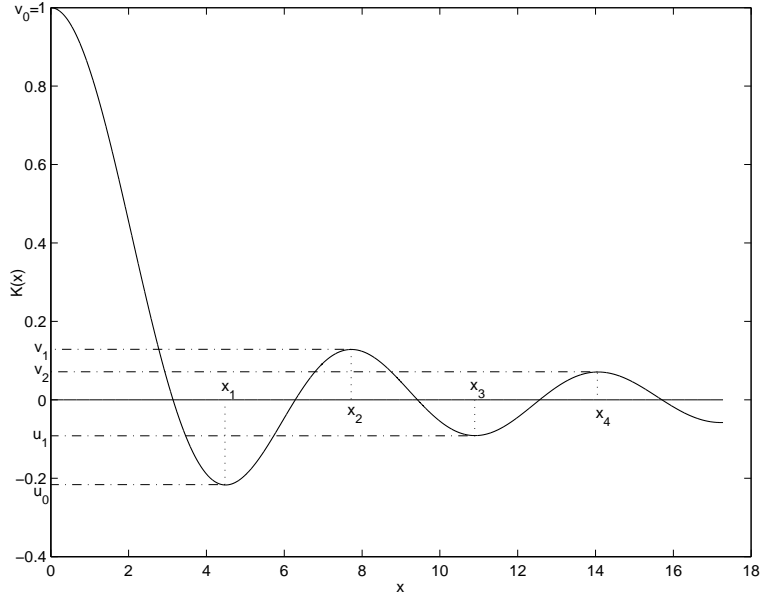}
\end{center}
\caption{The graphe of $K(x)$.}\label{functionKgraphe}
\end{figure}

Consider the sequence
\begin{equation}\label{sequencexk}
\{x_k\}_{k\in\mathbb{N}}:=\{x\geq0 \ ; \ x=\tan(x) \},
\end{equation}
with
\begin{displaymath}
0=x_0<x_1< \cdots <x_k<\cdots.
\end{displaymath}
In fact, one can check that
\begin{displaymath}
\{x_k\}_{k\in\mathbb{N}}=\{x\geq0 \ ; \ K^{\prime}(x)=0 \}.
\end{displaymath}
Moreover, for all $k\in\mathbb{N}^*$,
\begin{displaymath}
x_k\in(k\pi,k\pi+\frac{\pi}{2}).
\end{displaymath}
Define two sequences $\{u_k\}$ and $\{v_k\}$, for
$k\in\mathbb{N}$, by
\begin{displaymath}
u_k:=\cos(x_{2k+1})<0, \qquad v_k:=\cos(x_{2k})>0.
\end{displaymath}
Using the definition of $x_k$, one can see that
\begin{displaymath}
u_k=K(x_{2k+1}) \qquad \textrm{ and } \qquad v_k=K(x_{2k}).
\end{displaymath}
Thus, the sequence $\{u_k\}_{k\in\mathbb{N}}$ is increasing with
$-1<u_k<0$ and the sequence $\{v_k\}_{k\in\mathbb{N}}$ is
decreasing with $v_0=1$ and $0<v_k<1/2$ for $k\geq 1$ (see Figure
\ref{functionKgraphe}).

\noindent Moreover,
\begin{displaymath}
\lim_{k\to +\infty} u_k = \lim_{k\to +\infty} v_k = 0.
\end{displaymath}
Furthermore, one can check that, as soon as $\beta^*<0$,
\begin{displaymath}
\frac{\delta+\beta^*}{2\beta^*}<1=v_0.
\end{displaymath}
Finally, define a function $h$, for $x\in[-1,1/2)$, by
\begin{displaymath}
h(x)=\frac{2(1-x)}{1-2x}
\end{displaymath}
and set
\begin{displaymath}
h(v_0)=+\infty.
\end{displaymath}
We have the following results about the properties of the function
$h.$

%%
%%      LEMMA
%%
\begin{lemma}\label{lemma}
Suppose that
\begin{displaymath}
h(u_0)\leq n\frac{\beta_0-\delta}{\beta_0} \qquad \textrm{ and }
\qquad \delta+\beta^*\neq0.
\end{displaymath}
(i) If $\delta+\beta^*>0$, then there exists $k\in\mathbb{N}$ such
that
\begin{displaymath}
h(u_k)\leq n\frac{\beta_0-\delta}{\beta_0}<h(u_{k+1}).
\end{displaymath}
(ii) If $\delta+\beta^*<0$, then there exists $k\in\mathbb{N}$
such that
\begin{displaymath}
h(v_{k+1})\leq n\frac{\beta_0-\delta}{\beta_0}<h(v_k).
\end{displaymath}
\end{lemma}

\begin{proof}
Since the function $h$ is increasing on the interval $[-1,1/2)$,
we can see that
\begin{displaymath} h(u_k)\leq
n\frac{\beta_0-\delta}{\beta_0}<h(u_{k+1})
\end{displaymath}
is equivalent to
\begin{displaymath}
u_k\leq \frac{\delta+\beta^*}{2\beta^*}<u_{k+1}
\end{displaymath}
and
\begin{displaymath}
h(v_{k+1})\leq n\frac{\beta_0-\delta}{\beta_0}<h(v_k)
\end{displaymath}
is equivalent to
\begin{displaymath}
v_{k+1}\leq \frac{\delta+\beta^*}{2\beta^*}<v_k.
\end{displaymath}
The lemma now follows.
\end{proof}

%%
%%      PROPOSITION: EXISTENCE OF PURELY IMAGINARY ROOTS
%%
\begin{proposition}\label{propsolutionsystem}
(i) If
\begin{displaymath}
h(u_k)<n\frac{\beta_0-\delta}{\beta_0}<h(u_{k+1}),\qquad
k\in\mathbb{N},
\end{displaymath}
then system (\ref{eqsin})-(\ref{eqcos}) has exactly $2(k+1)$
solutions $(\tau_{1,1},\omega_{1,1})$, \dots,
$(\tau_{k+1,1},\omega_{k+1,1})$ and $(\tau_{1,2},\omega_{1,2})$,
\dots, $(\tau_{k+1,2},\omega_{k+1,2})$ with
\begin{displaymath}
\left\{ \begin{array}{ll}
\omega_{l,1}\tau_{l,1}\in((2l-1)\pi, x_{2l-1}), &\quad\textrm{ for }
l=1,\dots,k+1,\\
\omega_{l,2}\tau_{l,2}\in(x_{2l-1}, 2l\pi), &\quad\textrm{ for }
l=1,\dots,k+1
\end{array}\right.
\end{displaymath}
and
\begin{displaymath}
0<\tau_{1,1}<\cdots<\tau_{k+1,1}<\tau_{k+1,2}<\cdots<\tau_{1,2}.
\end{displaymath}
(ii) If
\begin{displaymath}
n\frac{\beta_0-\delta}{\beta_0}=h(u_k),\qquad k\in\mathbb{N},
\end{displaymath}
then system (\ref{eqsin})-(\ref{eqcos}) has exactly $2k+1$
solutions $(\tau_{1,1},\omega_{1,1})$, \dots,
$(\tau_{k+1,1},\omega_{k+1,1})$ and $(\tau_{1,2},\omega_{1,2})$,
\dots, $(\tau_{k,2},\omega_{k,2})$ with
\begin{displaymath}
\left\{ \begin{array}{ll}
\omega_{l,1}\tau_{l,1}\in((2l-1)\pi, x_{2l-1}),&\quad\textrm{ for } l=1,\dots,k,\\
\omega_{l,2}\tau_{l,2}\in(x_{2l-1}, 2l\pi),&\quad\textrm{ for } l=1,\dots,k,\\
\omega_{k+1,1}\tau_{k+1,1}=x_{2k+1}&
\end{array}\right.
\end{displaymath}
and
\begin{displaymath}
0<\tau_{1,1}<\cdots<\tau_{k+1,1}<\tau_{k,2}<\cdots<\tau_{1,2}.
\end{displaymath}
(iii) If
\begin{displaymath}
h(v_{k+1})<n\frac{\beta_0-\delta}{\beta_0}<h(v_{k}),\qquad
k\in\mathbb{N}^*,
\end{displaymath}
then system (\ref{eqsin})-(\ref{eqcos}) has exactly $2k+1$
solutions $(\tau_{1,1},\omega_{1,1})$, \dots,
$(\tau_{k+1,1},\omega_{k+1,1})$ and $(\tau_{1,2},\omega_{1,2})$,
\dots, $(\tau_{k,2},\omega_{k,2})$ with
\begin{displaymath}
\left\{ \begin{array}{ll}
\omega_{1,1}\tau_{1,1}\in(\pi/2, \pi),&\\
\omega_{l,1}\tau_{l,1}\in(x_{l+1}, (l+2)\pi),&\qquad \textrm{for } l=2,\dots,k+1,\\
\omega_{l,2}\tau_{l,2}\in((l+1)\pi,x_{l+1}),&\qquad \textrm{for }
l=1,\dots,k
\end{array}\right.
\end{displaymath}
and
\begin{displaymath}
0<\tau_{1,1}<\cdots<\tau_{k+1,1}<\tau_{k,2}<\cdots<\tau_{1,2}.
\end{displaymath}
(iv) If
\begin{displaymath}
n\frac{\beta_0-\delta}{\beta_0}=h(v_{k}),\qquad k\in\mathbb{N}^*,
\end{displaymath}
then system (\ref{eqsin})-(\ref{eqcos}) has exactly $2k$ solutions
$(\tau_{1,1},\omega_{1,1})$, \dots, $(\tau_{k,1},\omega_{k,1})$
and $(\tau_{1,2},\omega_{1,2})$, \dots,
$(\tau_{k,2},\omega_{k,2})$, with
\begin{displaymath}
\left\{ \begin{array}{ll}
\omega_{1,1}\tau_{1,1}\in(\pi/2, \pi),&\\
\omega_{l,1}\tau_{l,1}\in(x_l, (l+1)\pi),&\quad \textrm{for } l=2,\dots,k,\\
\omega_{l,2}\tau_{l,2}\in((l+1)\pi, x_{l+1}),&\quad \textrm{for }
l=1,\dots,k-1,\\
\omega_{k,2}\tau_{k,2}=x_{2k}&
\end{array}\right.
\end{displaymath}
and
\begin{displaymath}
0<\tau_{1,1}<\cdots<\tau_{k,1}<\tau_{k,2}<\cdots<\tau_{1,2}.
\end{displaymath}
(v) If
\begin{displaymath}
h(v_{1})<n\frac{\beta_0-\delta}{\beta_0}<h(v_{0}),
\end{displaymath}
then system (\ref{eqsin})-(\ref{eqcos}) has a unique solution
$(\tau_{1},\omega_{1})$ such that $\tau_1>0$ and
\begin{displaymath}
\omega_1\tau_1\in(0, \pi).
\end{displaymath}
\end{proposition}

\begin{proof}
We only prove $(i)$ when $k=0$. The other cases can be deduced
similarly. Assume that
\begin{displaymath}
h(u_0)<n\frac{\beta_0-\delta}{\beta_0}<h(u_1).
\end{displaymath}
This is equivalent to
\begin{displaymath}
u_0<\frac{\delta+\beta^*}{2\beta^*}<u_1.
\end{displaymath}
The function $K$ is strictly negative and decreasing on
$(\pi,x_1)$ with $K(y)\in(u_0,0)$ (see Figure
\ref{functionKgraphe}). So the equation
\begin{displaymath}
K(y)=\frac{\delta+\beta^*}{2\beta^*}
\end{displaymath}
has a unique solution $y_1$ on the interval $(\pi,x_1)$. Set
\begin{displaymath}
\tau_{1,1}=\frac{(y_1)^2}{2\beta^*(\cos(y_1)-1)}
\end{displaymath}
and $\omega_{1,1}=y_1/\tau_{1,1}$. Then,
$(\tau_{1,1},\omega_{1,1})$ is a unique solution of system
(\ref{eqsin})-(\ref{eqcos}) satisfying
$\omega_{1,1}\tau_{1,1}\in(\pi, x_1)$.

Moreover, the function $K$ is strictly negative and increasing on
$(x_1,2\pi)$ with $K(y)\in(u_0,0)$, so the equation
$K(y)=(\delta+\beta^*)/2\beta^*$ has a unique solution $y_2$ on
the interval $(x_1,2\pi)$. Set
\begin{displaymath}
\tau_{1,2}=\frac{(y_2)^2}{2\beta^*(\cos(y_2)-1)}
\end{displaymath}
and $\omega_{1,2}=y_2/\tau_{1,2}$. Then,
$(\tau_{1,2},\omega_{1,2})$ is a unique solution of system
(\ref{eqsin})-(\ref{eqcos}) which satisfies
$\omega_{1,2}\tau_{1,2}\in(x_1,2\pi)$.

Furthermore, the function $K$ is nonnegative on $[0,\pi]$ and
\begin{displaymath}
u_1=K(x_3)=\min_{x\geq2\pi}K(x).
\end{displaymath}
Therefore, system (\ref{eqsin})-(\ref{eqcos}) has two solutions,
$(\tau_{1,1},\omega_{1,1})$ and $(\tau_{1,2},\omega_{1,2})$.

Finally, using the fact that
\begin{displaymath}
\cos(y_1)\leq\cos(y_2),
\end{displaymath}
we obtain that \begin{displaymath} \tau_{1,1}<\tau_{1,2}.
\end{displaymath}
This completes the proof.
\end{proof}

Lemma \ref{lemma} and Proposition \ref{propsolutionsystem} give
conditions for the existence of pairs of purely imaginary roots of
equation (\ref{characteristicequationzero}). In the next
proposition, we study the properties of the purely imaginary roots
of (\ref{characteristicequationzero}).

%%
%%      PROPOSITION: PROPERTIES OF PURELY IMAGINARY ROOTS
%%
\begin{proposition}\label{prophopfbifurcation}
Assume that there exists a $\tau_c>0$ such that equation
(\ref{characteristicequationzero}) has a pair of purely imaginary
roots $\pm i\omega_c$ for $\tau=\tau_c$ with $\omega_c>0$. If
\begin{displaymath}
\omega_c\tau_c\neq x_k \qquad \textrm{ for all }k\in\mathbb{N},
\end{displaymath}
where the sequence $\{x_k\}_{k\in\mathbb{N}}$ is defined by
(\ref{sequencexk}), then $\pm i\omega_c$ are simple roots such
that
\begin{displaymath}
\left\{ \begin{array}{ll} \displaystyle\frac{d
Re(\lambda)}{d\tau}\Big|_{\tau=\tau_c}>0,& \qquad \textrm{ if }
\;\; \omega_c\tau_c\in(x_{2k},x_{2k+1}),\vspace{1ex}\\
\displaystyle \frac{d Re(\lambda)}{d\tau}\Big|_{\tau=\tau_c}<0,&
\qquad \textrm{ if } \;\; \omega_c\tau_c\in(x_{2k+1},x_{2k+2}),
\quad k\in\mathbb{N}.
\end{array}\right.
\end{displaymath}
\end{proposition}

\begin{proof}
Assume that there exists a $\tau_c>0$ such that equation
(\ref{characteristicequationzero}) has a pair of purely imaginary
roots $\pm i\omega_c$ for $\tau=\tau_c$ with $\omega_c>0$. Then,
$\omega_c\tau_c$ satisfies system (\ref{eqsin})-(\ref{eqcos}).

Assume that
\begin{displaymath}
\omega_c\tau_c\neq x_k \quad \textrm{ for all }k\in\mathbb{N}.
\end{displaymath}
Let us show that $\pm i\omega_c$ are simple characteristic roots
of (\ref{characteristicequationzero}). Using
(\ref{firstderivativezero}), one can see that $\pm i\omega_c$ are
simple roots of (\ref{characteristicequationzero}) if
\begin{displaymath}
1+2\beta^*\displaystyle\frac{\partial S}{\partial
\omega}(\tau_c,\omega_c)\neq 0 \qquad \textrm{ or } \qquad
\displaystyle\frac{\partial C}{\partial
\omega}(\tau_c,\omega_c)\neq 0.
\end{displaymath}
We will show that
\begin{displaymath}
\frac{\partial C}{\partial \omega}(\tau_c,\omega_c)\neq0.
\end{displaymath}
A simple computation shows that
\begin{displaymath}
\frac{\partial C}{\partial
\omega}(\tau_c,\omega_c)=\frac{g(\omega_c\tau_c)}{\omega_c^2\tau_c},
\end{displaymath}
where the function $g$ is defined by
\begin{displaymath}
g(x)=x\cos(x)-\sin(x) \qquad\textrm{for } x\geq0.
\end{displaymath}
One can check that $g(x)=0$ if and only if there exists a
$k_0\in\mathbb{N}$ such that $x=x_{k_0}$. Moreover,
\begin{displaymath}
g(x)>0 \qquad \textrm{ if and only if } \qquad
x\in(x_{2k+1},x_{2k+2}), \ k\in\mathbb{N}.
\end{displaymath}
This yields that
\begin{displaymath}
\frac{\partial C}{\partial \omega}(\tau_c,\omega_c)<0 \quad
\textrm{ if } \quad \omega_c\tau_c\in(x_{2k},x_{2k+1})
\end{displaymath}
and
\begin{displaymath}
\frac{\partial C}{\partial \omega}(\tau_c,\omega_c)>0 \quad
\textrm{ if } \quad \omega_c\tau_c\in(x_{2k+1},x_{2k+2}).
\end{displaymath}
Hence, $\pm i\omega_c$ are simple characteristic roots of
(\ref{characteristicequationzero}).

Let $\lambda(\tau)=\mu(\tau)+i\omega(\tau)$ be a characteristic
root of (\ref{characteristicequationzero}) such that
$\lambda(\tau_c)=\pm i\omega_c$. By separating the real and
imaginary parts, we obtain that
\begin{displaymath}
\left\{\begin{array}{rcl}
\mu(\tau)+\delta+\beta^*-\displaystyle\frac{2\beta^*}{\tau}
\int_0^{\tau}e^{-\mu(\tau) r}\cos(\omega(\tau) r)dr&=&0,\vspace{1ex}\\
\omega(\tau)+\displaystyle\frac{2\beta^*}{\tau}\int_0^{\tau}e^{-\mu(\tau)
r}\sin(\omega(\tau) r )dr&=&0.
\end{array}\right.
\end{displaymath}
We denote by $\mu^{\prime}(\tau)$ (respectively
$\omega^{\prime}(\tau)$) the first derivative of $\mu(\tau)$
(respectively $\omega(\tau)$) with respect to $\tau$. For
$\tau=\tau_c$, we obtain that
\begin{equation}\label{muprimezero}
\begin{array}{l}
\mu^{\prime}(\tau_c)\Big[1+2\beta^*\displaystyle\frac{\partial S}
{\partial \omega}(\tau_c,\omega_c)\Big]\vspace{1ex}\\
\qquad =2\beta^*\displaystyle\frac{\partial C}{\partial
\omega}(\tau_c,\omega_c)\omega^{\prime}(\tau_c)\vspace{1ex}
+\displaystyle\frac{2\beta^*}{\tau_c}\big(\cos(\omega_c\tau_c)-C(\tau_c,\omega_c)\big)
\end{array}
\end{equation}
and
\begin{equation}\label{omegaprimezero}
\begin{array}{l}
\omega^{\prime}(\tau_c)\Big[1+2\beta^*\displaystyle\frac{\partial S}{\partial \omega}(\tau_c,\omega_c)\Big]\vspace{1ex}\\
\qquad =-2\beta^*\displaystyle\frac{\partial C}{\partial
\omega}(\tau_c,\omega_c)\mu^{\prime}(\tau_c)\vspace{1ex}
+\displaystyle\frac{2\beta^*}{\tau_c}\big(S(\tau_c,\omega_c)-\sin(\omega_c\tau_c)\big).
\end{array}
\end{equation}
We consider two cases. First, assume that
\begin{equation}\label{hypS}
1+2\beta^*\displaystyle\frac{\partial S}{\partial
\omega}(\tau_c,\omega_c)=0.
\end{equation}
One can verify that
\begin{displaymath}
1+2\beta^*\displaystyle\frac{\partial S}{\partial
\omega}(\tau_c,\omega_c)=2+(\delta+\beta^*)\tau_c.
\end{displaymath}
Consequently, (\ref{hypS}) is equivalent to
\begin{equation}\label{tauc}
\tau_c=-\frac{2}{\delta+\beta^*}.
\end{equation}
Then, it follows from equation (\ref{omegaprimezero}) that
\begin{displaymath}
\frac{\partial C}{\partial
\omega}(\tau_c,\omega_c)\mu^{\prime}(\tau_c)=\frac{S(\tau_c,\omega_c)-\sin(\omega_c\tau_c)}{\tau_c}.
\end{displaymath}
Moreover, by using (\ref{eqsin}) and (\ref{eqcos}), we have
\begin{displaymath}
\begin{array}{rcl}
\displaystyle\frac{S(\tau_c,\omega_c)-\sin(\omega_c\tau_c)}{\tau_c}
&=&\displaystyle\frac{1-\big(\cos(\omega_c\tau_c)+\omega_c\tau_c\sin(\omega_c\tau_c)\big)}{\omega_c\tau_c^2}\\
&=&-\displaystyle\frac{\delta+\beta^*}{4\beta^*}\omega_c.
\end{array}
\end{displaymath}
Hence, (\ref{tauc}) implies that
\begin{displaymath}
\frac{\partial C}{\partial
\omega}(\tau_c,\omega_c)\mu^{\prime}(\tau_c)=\frac{\omega_c}{2\beta^*\tau_c}<0.
\end{displaymath}
Since $\frac{\partial C}{\partial \omega}(\tau_c,\omega_c)\neq0,$
we have
\begin{displaymath}
\mu^{\prime}(\tau_c)\neq0.
\end{displaymath}
Furthermore, the sign of $\mu^{\prime}(\tau_c)$ is the same as the
sign of $-\frac{\partial C}{\partial \omega}(\tau_c,\omega_c).$

We now assume that
\begin{displaymath}
1+2\beta^*\displaystyle\frac{\partial S}{\partial
\omega}(\tau_c,\omega_c)\neq0.
\end{displaymath}
Then, by using (\ref{muprimezero}) and (\ref{omegaprimezero}), we
obtain that $\mu^{\prime}(\tau_c)$ satisfies
\begin{displaymath}
\begin{array}{l}
\mu^{\prime}(\tau_c)\Big[\Big(1+2\beta^*\displaystyle\frac{\partial S}
{\partial \omega}(\tau_c,\omega_c)\Big)^2+\Big(2\beta^*\displaystyle\frac{\partial C}
{\partial \omega}(\tau_c,\omega_c)\Big)^2\Big]\vspace{1ex}\\
\quad =\displaystyle\frac{2\beta^*}{\tau_c}\Big[
2\beta^*\displaystyle\frac{\partial C}{\partial \omega}(\tau_c,\omega_c)
\big(S(\tau_c,\omega_c)-\sin(\omega_c\tau_c)\big)\vspace{1ex}\\
\qquad\qquad+\Big(1+2\beta^*\displaystyle\frac{\partial S}
{\partial \omega}(\tau_c,\omega_c)\Big)\big(\cos(\omega_c\tau_c)
-C(\tau_c,\omega_c)\big) \Big].
\end{array}
\end{displaymath}
Using the definitions of $C$ and $S$, one can check that
\begin{displaymath}
\frac{\partial C}{\partial
\omega}(\tau_c,\omega_c)\big(S(\tau_c,\omega_c)-\sin(\omega_c\tau_c)\big)+\frac{\partial
S}{\partial
\omega}(\tau_c,\omega_c)\big(\cos(\omega_c\tau_c)-C(\tau_c,\omega_c)\big)=0.
\end{displaymath}
Hence,
\begin{displaymath}
\mu^{\prime}(\tau_c)\Big[\Big(1+2\beta^*\displaystyle\frac{\partial
S}{\partial
\omega}(\tau_c,\omega_c)\Big)^2+\Big(2\beta^*\displaystyle\frac{\partial
C}{\partial
\omega}(\tau_c,\omega_c)\Big)^2\Big]=\displaystyle\frac{2\beta^*}{\tau_c}\big(\cos(\omega_c\tau_c)-C(\tau_c,\omega_c)\big).
\end{displaymath}
Notice that
\begin{displaymath}
\frac{\cos(\omega_c\tau_c)-C(\tau_c,\omega_c)}{\tau_c}=\frac{g(\omega_c\tau_c)}{\omega_c\tau_c^2}=\frac{\omega_c}{\tau_c}\frac{\partial
C}{\partial \omega}(\tau_c,\omega_c)\neq0.
\end{displaymath}
Since $1+2\beta^*\displaystyle\frac{\partial S}{\partial
\omega}(\tau_c,\omega_c)\neq0$, it follows that
\begin{displaymath}
\Big(1+2\beta^*\displaystyle\frac{\partial S}{\partial
\omega}(\tau_c,\omega_c)\Big)^2+\Big(2\beta^*\displaystyle\frac{\partial
C}{\partial \omega}(\tau_c,\omega_c)\Big)^2>0.
\end{displaymath}
Consequently, $\mu^{\prime}(\tau_c)\neq0$, and the sign of
$\mu^{\prime}(\tau_c)$ is the same as the sign of $-\frac{\partial
C}{\partial \omega}(\tau_c,\omega_c).$

In summary, we have obtained that, for $\tau=\tau_c$, equation
(\ref{characteristicequationzero}) has a pair of simple purely
imaginary roots $\pm i\omega_c$ such that
\begin{displaymath}
\left\{ \begin{array}{ll} \displaystyle\frac{d
Re(\lambda)}{d\tau}\Big|_{\tau=\tau_c}>0,& \qquad \textrm{ if }
\omega_c\tau_c\in(x_{2k},x_{2k+1}),\vspace{1ex}\\
\displaystyle \frac{d Re(\lambda)}{d\tau}\Big|_{\tau=\tau_c}<0,&
\qquad \textrm{ if } \omega_c\tau_c\in(x_{2k+1},x_{2k+2}).
\end{array}\right.
\end{displaymath}
This completes the proof.
\end{proof}

\begin{remark}
{\rm If there exists a $k\in\mathbb{N}^*$ such that
\begin{displaymath}
\omega_c\tau_c=x_k,
\end{displaymath}
then either $\pm i\omega_c$ are not simple roots of
(\ref{characteristicequationzero}) or
\begin{displaymath}
\frac{d Re(\lambda)}{d\tau}\Big|_{\tau=\tau_c}=0.
\end{displaymath}
Using a similar argument as in the proof of Proposition
\ref{propsolutionsystem}, we obtain that
\begin{displaymath}
\frac{\partial C}{\partial \omega}(\tau_c,\omega_c)=0.
\end{displaymath}
Thus, if
\begin{displaymath}
\tau_c=-\frac{2}{\delta+\beta^*},
\end{displaymath}
then $\pm i\omega_c$ are not simple roots of
(\ref{characteristicequationzero}). If
\begin{displaymath}
\tau_c\neq-\frac{2}{\delta+\beta^*},
\end{displaymath}
then
\begin{displaymath}
\frac{d Re(\lambda)}{d\tau}\Big|_{\tau=\tau_c}=0.
\end{displaymath}}
\end{remark}

In the next theorem, we show that there exists a Hopf bifurcation
at the non-trivial equilibrium $x\equiv x^*$ of equation
(\ref{equationx}).

%%
%%      THEOREM: HOPF BIFURCATIONS
%%
\begin{theorem}\label{theoremhopfbifurcation}
Assume that
\begin{displaymath}
h(u_0)\leq n\frac{\beta_0-\delta}{\beta_0} \qquad \textrm{ and }
\qquad \delta+\beta^*\neq0.
\end{displaymath}
Then a Hopf bifurcation occurs at $x\equiv x^*$ for
$\tau=\tau_0:=\min_{\omega_c\tau_c\neq x_k, \ k\in\mathbb{N}}
\tau_c$, where $(\tau_c,\omega_c)$ are solutions of
(\ref{eqsin})-(\ref{eqcos}), defined in Proposition
\ref{propsolutionsystem}. When $0\leq\tau<\tau_0$, the equilibrium
$x\equiv x^*$ is locally asymptotically stable and it is unstable
while $\tau_0\leq\tau\leq\tau_l$, where $\tau_l$ is the larger
value of $\tau_c$ such that $\omega_c\tau_c\in(x_{2k},x_{2k+1})$,
$k\in\mathbb{N}$.
\end{theorem}

\begin{proof}
We first check that $x\equiv x^*$ is locally asymptotically stable
when $\tau\in[0,\tau_0)$. Notice that when $\tau\in[0,\tau_0)$,
equation (\ref{characteristicequationzero}) does not have purely
imaginary roots. Let $\tau^*>0$ be small enough and fixed. Assume
that, for $\tau\in(0,\tau^*)$, equation
(\ref{characteristicequationzero}) has a characteristic root
$\lambda(\tau)=\mu(\tau)+i\omega(\tau)$ with $\mu(\tau)>0$.
Separating the real and imaginary parts, we obtain
\begin{displaymath}
\mu(\tau)=-(\delta+\beta^*)+\frac{2\beta^*}{\tau}\int_0^{\tau}e^{-\mu(\tau)r}\cos(\omega(\tau)
r)dr
\end{displaymath}
and
\begin{displaymath}
\omega(\tau)=-\frac{2\beta^*}{\tau}\int_0^{\tau}e^{-\mu(\tau)r}\sin(\omega(\tau)
r)dr.
\end{displaymath}
We deduce that, for $\tau\in(0,\tau^*)$,
\begin{displaymath}
|\mu(\tau)|\leq |\delta+\beta^*|-2\beta^* \qquad \textrm{ and }
\qquad |\omega(\tau)|\leq -2\beta^*.
\end{displaymath}
Consequently,
\begin{displaymath}
\lim_{\tau\to0} \tau\mu(\tau)=0 \qquad \textrm{ and } \qquad
\lim_{\tau\to0} \tau\omega(\tau)=0.
\end{displaymath}
Integrating by parts, we obtain
\begin{displaymath}
2\mu(\tau)=-(\delta+\beta^*)+2\beta^*e^{-\tau\mu(\tau)}K(\tau\omega(\tau)).
\end{displaymath}
Since $\mu(\tau)>0$, we have for $\tau\in(0,\tau^*)$ that
\begin{displaymath}
-(\delta+\beta^*)+2\beta^*e^{-\tau\mu(\tau)}K(\tau\omega(\tau))>0.
\end{displaymath}
When $\tau$ tends to zero, we obtain
\begin{displaymath}
\beta^*-\delta\geq0.
\end{displaymath}
However, $\beta^*-\delta<0$. This is a contradiction. Therefore,
for $\tau\in(0,\tau^*)$, $\mu(\tau)<0$. Applying Rouch\'e's
Theorem [\ref{rouche}, p.248], we obtain that all characteristic
roots of (\ref{characteristicequationzero}) have negative real
parts when $\tau\in[0,\min(\tau_c))$. Therefore, $x\equiv x^*$ is
locally asymptotically stable.

Using Lemma \ref{lemma}, Propositions \ref{propsolutionsystem} and
\ref{prophopfbifurcation}, we conclude to the existence of
$\tau_l$. This concludes the proof.
\end{proof}

We illustrate the results of Theorem \ref{theoremhopfbifurcation}
in the next corollary.

\begin{figure}[!hpt]
\begin{center}\includegraphics[width=10cm,height=8cm]{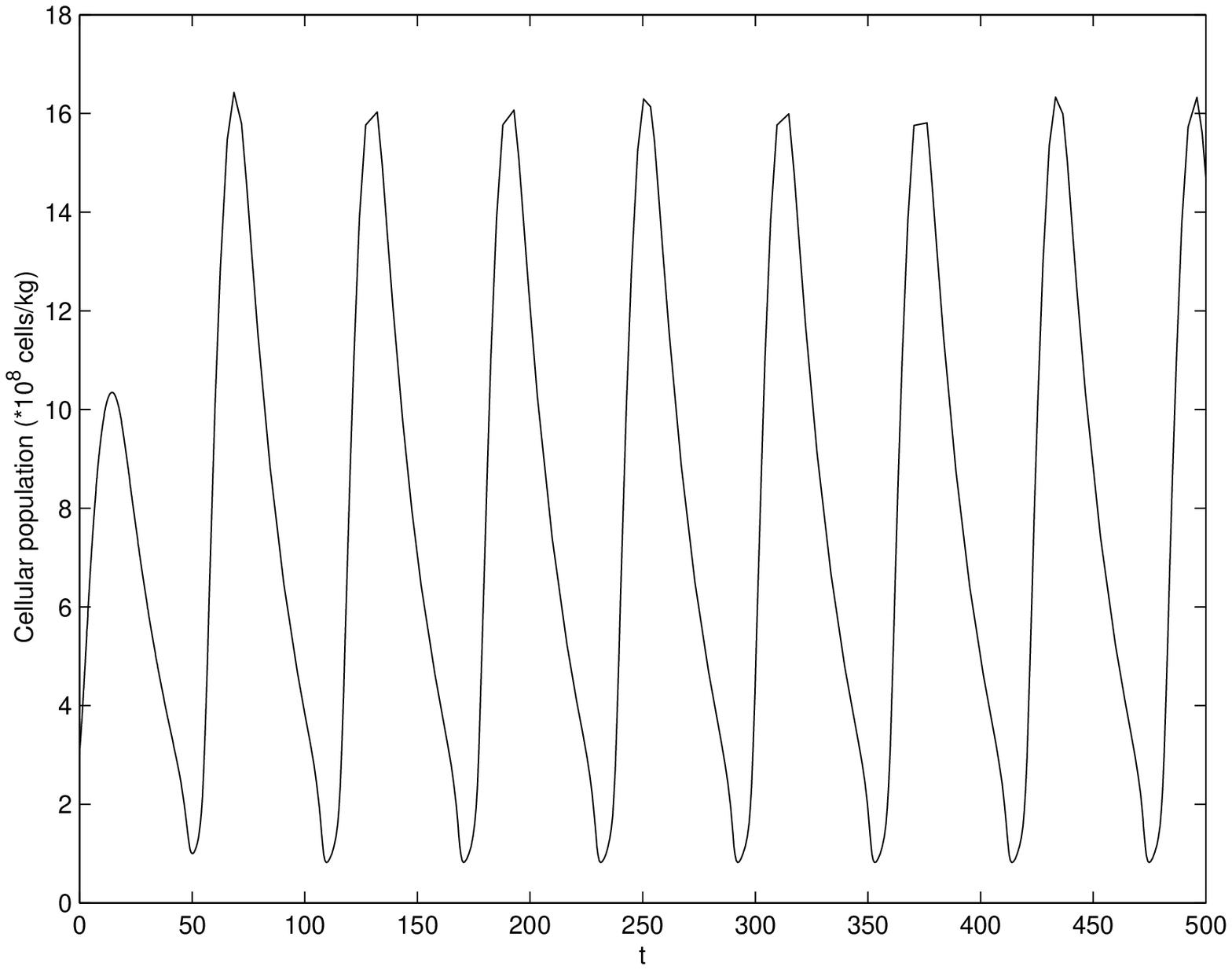}\end{center}
\begin{center}\includegraphics[width=10cm,height=8cm]{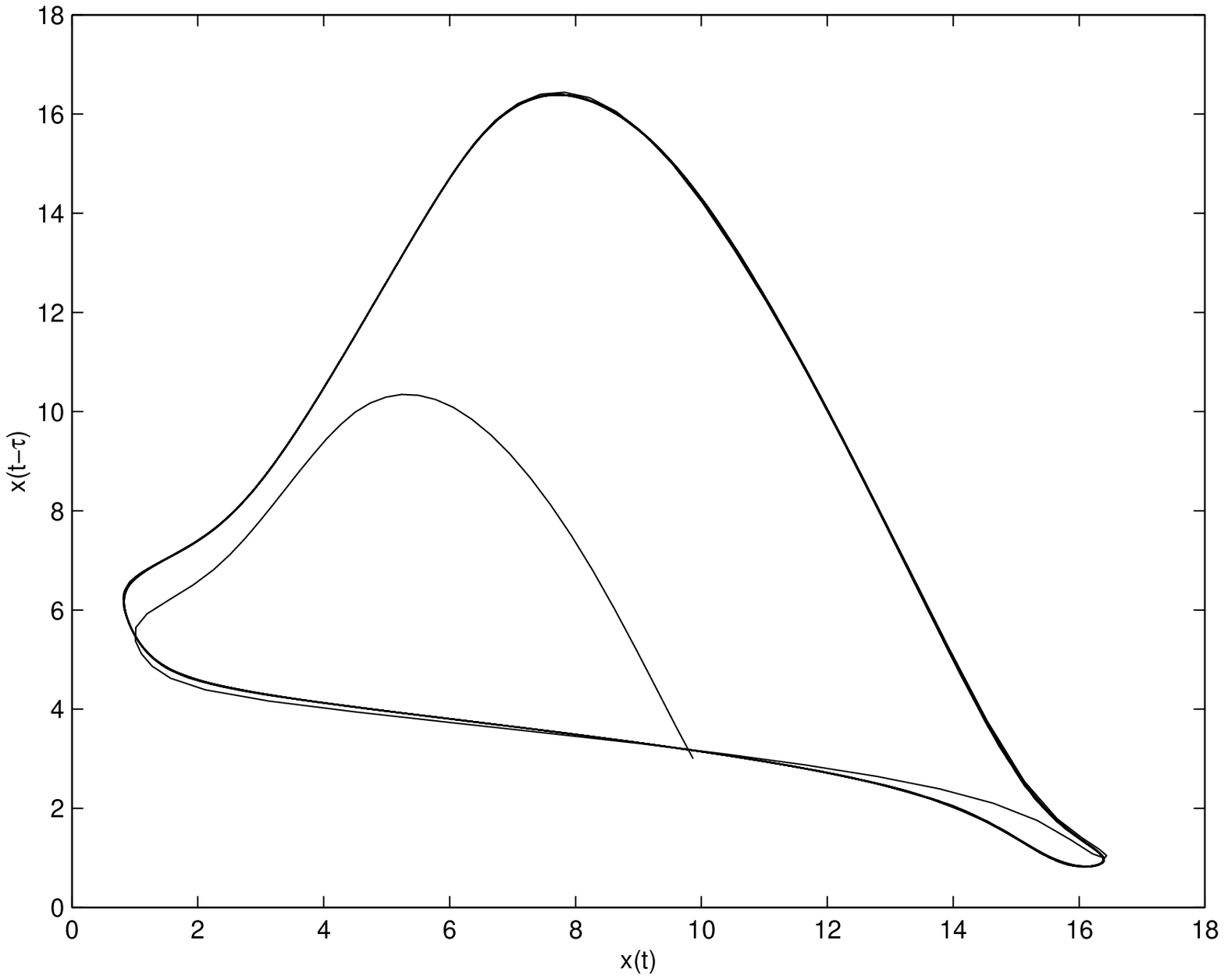}\end{center}
\caption{With the values given by (\ref{parametersvalues}) and
$\theta=1.62\times10^8$ cells/kg, equation (\ref{equationx}) has a
periodic solution for $\tau=18.2$ days. This solution has a period
about $50$ days. One can see that the solution reaches a limit
cycle.}\label{figbif}
\end{figure}

%%
%%      COROLLARY: APPLICATION
%%
\begin{corollary}
Assume that the parameters $\delta$, $\beta_0$ and $n$ are given
by (\ref{parametersvalues}). Then there exists a unique value
$\tau_c>0$ such that a Hopf bifurcation occurs at $x\equiv x^*$
when $\tau=\tau_c$. When $\tau<\tau_c$, the equilibrium is locally
asymptotically stable and becomes unstable when $\tau\geq \tau_c$.
Moreover, when $\tau=\tau_c$, equation (\ref{equationx}) has a
periodic solution with a period close to $46$ days (see Figure
\ref{figbif}). The value of $\tau_c$ is approximately given by
\begin{displaymath}
\tau_c\simeq 18 \textrm{ days}.
\end{displaymath}
\end{corollary}

\begin{proof}
With the values given by (\ref{parametersvalues}), we obtain
\begin{displaymath}
n\frac{\beta_0-\delta}{\beta_0}\simeq2.9153>h(v_1)\simeq 2.3455.
\end{displaymath}
Hence, Proposition \ref{propsolutionsystem} implies that the
system (\ref{eqsin})-(\ref{eqcos}) has a unique solution
$(\tau_c,\omega_c)$ with $\tau_c>0$ and
$\omega_c\tau_c\in(0,\pi)$. From Theorem
\ref{theoremhopfbifurcation}, we know that a Hopf bifurcation
occurs at $x\equiv x^*$ for $\tau=\tau_c$. The equilibrium is
locally asymptotically stable when $\tau<\tau_c$ and becomes
unstable when $\tau\geq \tau_c$. Consequently, for $\tau=\tau_c$,
equation (\ref{equationx}) has a periodic solution with a period
close to $2\pi/\omega_c$. One can check that
\begin{displaymath}
\tau_c\simeq 18 \textrm{ days } \qquad \textrm{ and } \qquad
\omega_c\simeq 0.138.
\end{displaymath}
Computer simulations confirm our analysis (see Figure
\ref{figbif}).
\end{proof}

As mentioned earlier, the results in Theorem
\ref{theoremhopfbifurcation} still hold when $\tau_{min}>0$.
However, in this case, the computations in the proof of Theorem
\ref{theoremhopfbifurcation} are much more complicated.

\section{Discussion}\label{discussion}

Hematological diseases have attracted a significant amount of
modeling attention because a number of them are periodic in nature
(Haurie {\it et al.} \cite{hdm98}). Some of these diseases involve
only one blood cell type and are due to the destabilization of
peripheral control mechanisms, e.g., periodic auto-immune
hemolytic anemia (B\'elair {\it et al.} \cite{bmm95} and Mahaffy
{\it et al.} \cite{mbm98}). Such periodic hematological diseases
involve periods between two and four times the bone marrow
production/maturation delay. Other periodic hematological
diseases, such as cyclical neutropenia (Haurie {\it et al.}
\cite{hdm98}), involve oscillations in all of the blood cells and
very long period dynamics on the order of weeks to months (Fowler
and Mackey \cite{fm02} and Pujo-Menjouet et al.
\cite{mackeypujo2}) and are thought to be due to a destabilization
of the pluripotent stem cell compartment from which all types of
mature blood cells are derived.

We have studied a scalar delay model that describes the dynamics
of a pluripotent stem cell population involved in the blood
production process in the bone marrow. The distributed delay
describes the cell cycle duration. We established stability
conditions for the model independent of the delay. We have also
observed oscillations in the pluripotent stem cell population
through Hopf bifurcations. With parameter values given in Mackey
\cite{mackey1978, mackey1979}, our calculations indicate that the
oscillatory pluripotent stem cell population involves a period of
$45$ days.

It will be very interesting to study the dynamics of the two
dimensional systems (Mackey \cite{mackey1978, mackey1979}, Mackey
{\it et al.} \cite{mackey2003}, Pujo-Menjouet {\it et al.}
\cite{mackeypujo2}) modeling the proliferating phase cells and
resting phase cells with distributed delays. We leave this for
future consideration.

\end{document}